\documentclass[11pt, a4paper]{article}

\usepackage{amsmath}
\usepackage{amssymb}
\usepackage{amsthm}
\usepackage{appendix}

\newtheorem{theorem}{Theorem}[section]

\newtheorem{remark}{Remark}[section]
\newtheorem{example}{Example}[section]

\begin{document}
\title{On exact systems $\{t^{\alpha}\cdot e^{2\pi i nt}\}_{n\in\mathbb{Z}\setminus A}$ in $L^2 (0,1)$\\ which are not Schauder Bases
and their generalizations}
\author{Elias Zikkos\\
Khalifa University of Science and Technology\\
Abu Dhabi, United Arab Emirates\\
email address:  elias.zikkos@ku.ac.ae and eliaszikkos@yahoo.com}

\maketitle

\begin{abstract}
Let $\{e^{i\lambda_n t}\}_{n\in\mathbb{Z}}$ be an exponential Schauder Basis for $L^2 (0,1)$, for $\lambda_n\in\mathbb{R}$,
and let $\{r_n(t)\}_{n\in\mathbb{Z}}$ be its dual Schauder Basis.
Let $A$ be a non-empty subset of the integers containing exactly $M$ elements. We prove that for $\alpha >0$ the weighted system
\[
\{t^{\alpha}\cdot r_n(t)\}_{n\in\mathbb{Z}\setminus A}
\]
is exact in the space $L^2 (0,1)$, that is, it is complete and minimal in $L^2 (0,1)$, if and only if
\[
M-\frac{1}{2}\le \alpha< M+\frac{1}{2}.
\]
We also show that such a system is not a Riesz Basis for $L^2 (0,1)$.
In particular, the weighted trigonometric system $\{t^{\alpha}\cdot e^{2\pi i n t}\}_{n\in\mathbb{Z}\setminus A}$ is exact
in $L^2 (0,1)$, if and only if
$\alpha\in [M-\frac{1}{2}, M+\frac{1}{2})$, but it is not a Schauder Basis for $L^2 (0,1)$.
\end{abstract}



\section{Introduction and our Result}
\setcounter{equation}{0}

In \cite{HeilYoon} Heil and Yoon investigated the properties of systems of weighted exponentials of the form
$\{g(t)\cdot e^{2\pi i n t}\}_{n\in\mathbb{Z}\setminus A}$ where $g\in L^2 (0,1)\setminus\{0\}$ and
$A$ is a non-empty finite subset of the integers containing exactly $M$ elements.
Amongst their results, they prove in \cite[Theorem 4]{HeilYoon} that the system
\[
\{t^{J}\cdot e^{2\pi i n t}\}_{n\in\mathbb{Z}\setminus A},
\]
where $J$ is a positive integer, is complete and minimal in $L^2 (0,1)$ if and only if $J=M$.

They also raise the interesting question \cite[page 100]{HeilYoon}
whether similar results hold for the system
\begin{equation}\label{system}
\{t^{\alpha}\cdot e^{2\pi i n t}\}_{n\in\mathbb{Z}\setminus A}
\end{equation}
where $\alpha$ is a real positive number, not necessarily an integer.

Motivated by the above and the work of Shukurov \cite{Shukurov}, we provide the following answer:
we prove that when the cardinality $|A|$ of the set $A$ is equal to $M$,
then the system $(\ref{system})$ is complete and minimal in $L^2 (0,1)$ if and only if
\begin{equation}\label{inequality}
M-\frac{1}{2}\le \alpha< M+\frac{1}{2}.
\end{equation}
We also show that such a complete and minimal system is not a Schauder Basis for $L^2 (0,1)$.

\begin{example}
The system $\{t^{\alpha}\cdot e^{2\pi i n t}\}_{n\in\mathbb{Z}\setminus A}$ is complete and minimal in $L^2 (0,1)$ if
\[
|A|=1,\qquad \text{and}\qquad \frac{1}{2}\le \alpha< \frac{3}{2},
\]
\[
|A|=2,\qquad \text{and}\qquad \frac{3}{2}\le \alpha< \frac{5}{2},
\]
\[
|A|=3,\qquad \text{and}\qquad \frac{5}{2}\le \alpha< \frac{7}{2},
\]
and so on, but no such system is a Schauder Basis for $L^2 (0,1)$.
\end{example}

We point out that the condition $(\ref{inequality})$ is necessary and sufficient in order for
more general systems of the form $\{t^{\alpha}\cdot r_n(t)\}_{n\in\mathbb{Z}\setminus A}$ to be complete and minimal in $L^2 (0,1)$,
where $\{r_n\}_{n\in\mathbb{Z}}$ is the dual Schauder basis for an exponential Schauder basis
$\{e^{i\lambda_n t}\}_{n\in\mathbb{Z}}$ for $L^2 (0,1)$. That is the content of our result which reads as follows.

\begin{theorem}\label{theorem}
Let $\{e^{i\lambda_n t}\}_{n\in\mathbb{Z}}$ be an exponential Schauder Basis for $L^2 (0,1)$, for $\lambda_n\in\mathbb{R}$,
and let $\{r_n\}_{n\in\mathbb{Z}}$ be its Dual Schauder Basis for $L^2 (0,1)$.
Let $A\subset \mathbb{Z}$ so that $|A|=M$ for some positive integer $M$. Then the system
\begin{equation}\label{system1}
\{t^{\alpha}\cdot r_n(t)\}_{n\in\mathbb{Z}\setminus A}
\end{equation}
is complete and minimal in the space $L^2 (0,1)$, if and only if $\alpha\in [M-\frac{1}{2}, M+\frac{1}{2})$.

However, for such values of $\alpha$ the system $(\ref{system1})$ is not a Riesz Basis for $L^2 (0,1)$ and neither a frame.

In particular, the weighted trigonometric system $\{t^{\alpha}\cdot e^{2\pi i n t}\}_{n\in\mathbb{Z}\setminus A}$
is complete and minimal in the space $L^2 (0,1)$, if and only if $\alpha\in [M-\frac{1}{2}, M+\frac{1}{2})$, but it
is not even a Schauder Basis for $L^2 (0,1)$.
\end{theorem}

We point out that under the assumption that the elements of the set $A$ are consecutive integers,
the non-basicity of the system $(\ref{system})$ follows from the work of Kazarian \cite[Corollary 1]{Kazarian1978}.
In such a case, a complete and minimal system $(\ref{system})$ is a Markushevich Basis \cite[Theorem 1]{Kazarian1987}.

\begin{remark}
The proof of Theorem $\ref{theorem}$ depends on the invertibility of Vandermonde matrices.
\end{remark}

\section{Preliminaries}
\setcounter{equation}{0}

For the terminology presented below, one may consult the books of Young \cite{Young},  Heil \cite{Heil}, and Christensen \cite{Christensen}.

Let $\mathcal{H}$ be a separable Hilbert space endowed with an inner product $\langle .\,  _,\, .\rangle$ and a
norm $||\, .\,  ||$. Let $F:=\{f_n\}_{n\in J}$ be a countable family of vectors in $\mathcal{H}$ with $J\subset \mathbb{Z}$.
We say that $F$ is complete if the closed span of $F$ in $\mathcal{H}$ is equal to $\mathcal{H}$ and we say that
$F$ is minimal if each element $f_n$ does not belong to the closed span of the remaining
vectors of $F$ in $\mathcal{H}$. If $F$ is both complete and minimal then it is called $\bf{exact}$. It is known that
the minimality of $F$ is equivalent to the existence of a biorthogonal sequence $\{g_n\}_{n\in J}$ in $\mathcal{H}$, that is
\[
\langle f_n , g_m \rangle = \begin{cases} 1, & n=m\\  0, & n\not=m \end{cases}.
\]
This biorthogonal family $\{g_n\}_{n\in J}$ is unique if $F$ is exact, but $\{g_n\}_{n\in J}$ itself does not have to be exact.
If both $F$ and $\{g_n\}_{n\in J}$ are exact, we say that they are Markushevich bases for $\mathcal{H}$.

An exact system $F$ is a Schauder Basis for $\mathcal{H}$, if for any $f\in \mathcal{H}$,
there exist unique numbers $d_{f,n}$ so that
\[
f=\sum_{n\in J} d_{f,n}  f_n,
\]
with convergence in the $L^2 (0,1)$ norm. In this case, the unique biorthogonal family $\{g_n\}_{n\in J}$ is also a Schauder Basis for
$\mathcal{H}$, called the dual Schauder basis of $F$, and for any $f\in \mathcal{H}$ we have
\[
f=\sum_{n\in J} \langle f, g_n \rangle f_n=\sum_{n\in J} \langle f, f_n \rangle g_n,
\]
with convergence in the $L^2 (0,1)$ norm (\cite[Theorem 5.12 and Corollary 5.22]{Heil} and \cite[relation (1)]{HeilYoon}).
A Schauder basis may converge conditionally,
whereas a Riesz Basis for $\mathcal{H}$ is a bounded Schauder basis that converges unconditionally (\cite[Theorem 7.11]{Heil}).
We recall that
$\{f_n\}_{n\in J}$  is a Riesz basis for $\mathcal{H}$ if $f_n=V(e_n)$ for all $n\in J$ where
$\{e_n\}_{n\in J}$ is an orthonormal basis for $\mathcal{H}$ and $V$ is a bounded
bijective operator from $\mathcal{H}$ onto $\mathcal{H}$.

\section{Proof of Theorem $\ref{theorem}$}
\setcounter{equation}{0}

First we will show that the system $(\ref{system1})$ is complete in $L^2 (0,1)$ if $\alpha \ge M - \frac{1}{2}$
and then show that it is minimal if $\alpha < M + \frac{1}{2}$. We then prove that for $\alpha\in [M - \frac{1}{2}, M + \frac{1}{2})$
the exact system $(\ref{system1})$ is not a Riesz basis for $L^2 (0,1)$, while the exact system $(\ref{system})$
is not even a Schauder Basis for $L^2 (0,1)$.

\subsection{If $\alpha\ge M-\frac{1}{2}$, the system $(\ref{system1})$ is Complete}

Since $\{r_n(t)\}_{n\in\mathbb{Z}}$ is the dual Schauder Basis of $\{e^{i\lambda_n t}\}_{n\in\mathbb{Z}}$ for $L^2 (0,1)$,
then for any $f\in L^2 (0,1)$ we have
\[
f(t)=\sum_{n\in\mathbb{Z}}\langle f, r_n\rangle e^{i\lambda_n t}
\]
with the series converging in the $L^2 (0,1)$ norm.

For notation purposes, let $t_{\alpha}(t):=t^{\alpha}$. Now,
if the system $(\ref{system1})$ is not complete in $L^2 (0,1)$, then there exists a non-zero function $h(t)\in L^2 (0,1)$ so that
\[
\langle h, t_{\alpha} r_n\rangle =0,\qquad \text{for\,\, all}\quad n\in \mathbb{Z}\setminus A.
\]
Since $|t^{\alpha}|\le 1$ for $t\in [0,1]$, then $ht_{\alpha}$ belongs to $L^2 (0,1)$. Combining the above shows that
\[
h(t)\cdot t^{\alpha}=\sum_{n\in A}\langle h t_{\alpha}, r_n\rangle e^{i\lambda_n t},
\]
almost everywhere on $[0,1]$, and such that at least one of the coefficients $\langle h t_{\alpha}, r_n\rangle$
for $n\in A$ is not equal to zero.

Rewrite now the set $\{\lambda_i:\,\, i\in A\}$ as $\{\Lambda_1, \Lambda_2, \dots ,\Lambda_M\}$
such that $\Lambda_1<\Lambda_2<\dots <\Lambda_M$. Then
\[
h(t)\cdot t^{\alpha}=\sum_{j=1}^{M} a_j e^{i\Lambda_j t},\qquad a_j\in\mathbb{C}.
\]
The analytic extension to the complex plane of the above right hand side is the entire function
\begin{equation}\label{trigpolc}
H_{\alpha}(z):=\sum_{j=1}^{M} a_j e^{i\Lambda_j z},
\end{equation}
so $H_{\alpha}(t)=h(t)\cdot t^{\alpha}$ for $t\in [0,1]$.
We will prove below that $H_{\alpha}$ cannot vanish $M$ times or more at the point $z=0$. Suppose it does: then
\begin{equation}\label{zeros}
0=H_{\alpha}(0)=H_{\alpha}^{'}(0)=H_{\alpha}^{''}(0)=\cdots H_{\alpha}^{(M-1)}(0)
\end{equation}
where $H_{\alpha}^{(k)}(z)$ is the $k^{th}$ derivative of $H_{\alpha}$ with respect to the variable $z$.
From $(\ref{trigpolc})$ and $(\ref{zeros})$, we get
\begin{eqnarray*}
0 & = & a_1 + a_2  +  a_3    +  a_4    +  \cdots  +  a_M\\
  & = & a_1\cdot \Lambda_1 + a_2\cdot \Lambda_2   +  a_3\cdot \Lambda_3   +  a_4\cdot \Lambda_4    +  \cdots  +  a_M\cdot \Lambda_M\\
  & = & a_1\cdot \Lambda_1^2 + a_2\cdot \Lambda_2^2 +  a_3\cdot \Lambda_3^2 +  a_4\cdot \Lambda_4^2 + \cdots  +  a_M\cdot \Lambda_M^2
\end{eqnarray*}
and so on, until the last equation ($0=H_{\alpha}^{(M-1)}(0)$)
\[
0 = a_1\cdot \Lambda_1^{M-1} + a_2\cdot \Lambda_2^{M-1} + a_3\cdot \Lambda_3^{M-1} + a_4\cdot \Lambda_4^{M-1} + \cdots  +  a_M\cdot \Lambda_M^{M-1}.
\]
We can write this system of equations as

\[
\begin{pmatrix}
  0 \\
  0  \\
  0 \\
 \vdots \\
 0
 \end{pmatrix}
 =
 \begin{pmatrix}
 1 & 1 & 1 & 1 & \cdots & 1 \\
 \Lambda_1  & \Lambda_2   & \Lambda_3 & \Lambda_4 &  \cdots & \Lambda_M  \\
 \Lambda_1^2 &  \Lambda_2^2 & \Lambda_3^2 & \Lambda_4^2 & \cdots & \Lambda_M^2  \\
  \Lambda_1^3 &  \Lambda_2^3 & \Lambda_3^3 & \Lambda_4^3 & \cdots & \Lambda_M^3  \\
 \vdots  & \vdots  & \vdots & \ddots & \vdots \\
\Lambda_1^{M-1} &  \Lambda_2^{M-1} & \Lambda_3^{M-1} & \Lambda_4^{M-1} & \cdots & \Lambda_M^{M-1}
 \end{pmatrix}
 \cdot
\begin{pmatrix}
  a_1\\
  a_2 \\
  a_3 \\
  a_4\\
  \vdots \\
  a_M
 \end{pmatrix}.
 \]
The above square matrix is an invertible Vandermonde matrix since the $\Lambda_i$'s are all different.
Clearly now we have $a_j=0$ for all $j=1, 2,3,\dots, M$, thus $H_{\alpha}(z)$ is the Zero function. But
$H_{\alpha}(t)=h(t)\cdot t^{\alpha}$ for $t\in [0,1]$, thus
$h(t)$ is identically equal to zero, a contradiction.
Therefore $H_{\alpha}(z)$ vanishes at the point $z=0$ at most $M-1$ times, hence
$H_{\alpha}(z)=z^k\cdot g(z)$ for some integer $k$ in the set $\{0,1,\dots, M-1\}$ and some entire function $g$ such that $g(0)\not =0$.
Since $H_{\alpha}(t)=h(t)\cdot t^{\alpha}$ for $t\in [0,1]$, we have
\[
h(t)\cdot t^{\alpha}=t^k\cdot g(t),\qquad \text{for\,\, some}\quad k\in\{0,1,2,\cdots, M-1\},\quad g\in C[0,1]\quad g(0)\not= 0.
\]
Thus
\[
|h(t)|^2=|g(t)|^2\cdot t^{2k-2\alpha}.
\]
Since $g(0)\not= 0$, then there are $\epsilon>0$ and $\delta>0$, so that $|g(t)|>\epsilon$ for all $t\in [0,\delta]$.
So, for any $0<\rho<\delta$, one has $|h(t)|^2>\epsilon^2\cdot t^{2k-2\alpha}$ for all $t\in [\rho ,\delta]$.
Therefore,
\[
\int_{\rho}^{\delta}|h(t)|^2\,\, dt> \epsilon^2\cdot \int_{\rho}^{\delta} t^{2k-2\alpha}\,\, dt.
\]
Since $k\le M-1$ and $\alpha\ge M-\frac{1}{2}$, then $2k-2\alpha\le -1$ and it follows that
\[
\int_{\rho}^{\delta} t^{2k-2\alpha}\,\, dt\to \infty\qquad\text{as}\qquad \rho\to 0^+.
\]
Hence $\int_{0}^{1}|h(t)|^2\,\, dt=\infty$ which contradicts the assumption that $h\in L^2 (0,1)$.
Therefore, if $\alpha\ge M-\frac{1}{2}$ then the system $(\ref{system1})$ is complete in $L^2(0,1)$.

\subsection{If $\alpha<M+\frac{1}{2}$, the system $(\ref{system1})$ is Minimal}

Like before, the set $\{\lambda_i:\,\, i\in A\}$ is rewritten as $\{\Lambda_1, \Lambda_2, \dots ,\Lambda_M\}$
such that $\Lambda_1<\Lambda_2<\dots <\Lambda_M$.

We will show that for each $n\in\mathbb{Z}\setminus A$, there exist real numbers $a_{n,j}$ for $j=1,2,\cdots, M$,
so that the entire function
\begin{equation}\label{trigpol2}
f_n(z):=e^{i\lambda_n z}+a_{n,1}e^{i\Lambda_1 z}+a_{n,2}e^{i\Lambda_2 z}+\cdots +a_{n,M}e^{i\Lambda_M z}
\end{equation}
vanishes at the point $z=0$ at least $M$ times and

\begin{equation}\label{Mc}
f_n(t)=t^M\cdot g_n(t)\qquad \text{for\,\, some\,\, continuous\,\, function}\quad  g_n\,\, \text{on}\,\, [0,1].
\end{equation}
The minimality of $(\ref{system1})$ follows from $(\ref{trigpol2})$ and $(\ref{Mc})$. Indeed, from $(\ref{Mc})$ we have
\[
\frac{f_n(t)}{t^{\alpha}}=
\frac{f_n(t)}{t^M}\cdot \frac{t^M}{t^{\alpha}}=g_n(t)\cdot \frac{t^M}{t^{\alpha}}.
\]
So
\[
\left|\frac{f_n(t)}{t^{\alpha}}\right|^2=|g_n(t)|^2\cdot t^{2M-2\alpha}.
\]
Since $g_n$ is continuous on $[0,1]$ and $\alpha<M+\frac{1}{2}$, thus $2M-2\alpha>-1$,
it follows that $f_n(t)/t^{\alpha}$ belongs to $L^2 (0,1)$. Then, for all positive integers $n, \, m$ which are not in the set $A$, we have
\[
\langle \frac{f_n}{t_{\alpha}}, t_{\alpha} r_m  \rangle = \int_{0}^{1}
\left(e^{i\lambda_n t}+a_{n,1}e^{i\Lambda_1 t}+a_{n,2}e^{i\Lambda_2 t}+\cdots +a_{n,M}e^{i\Lambda_M t}\right)\cdot\overline{r_m}\,\, dt=
\begin{cases} 1, & n=m\\  0, & n\not=m \end{cases}.
\]
Thus, the family $\{\frac{f_n(t)}{t^{\alpha}}\}_{n\in\mathbb{Z}\setminus A}$
is biorthogonal to the sequence $\{t^{\alpha}\cdot r_n(t)\}_{n\in\mathbb{Z}\setminus A}$.
This means that the system $(\ref{system1})$ is minimal.

\smallskip

So, it remains to verify the existence of the numbers $a_{n,j}$ so that $(\ref{Mc})$ is true. Now,
if a trigonometric polynomial $f_n$ as in $(\ref{trigpol2})$ vanishes at the point $z=0$ at least $M$ times, then
\begin{equation}\label{zeros1}
0=f_n(0)=f_n^{'}(0)=f_n^{''}(0)=\cdots f_n^{(M-1)}(0)
\end{equation}
where $f_n^{(k)}(z)$ is the $k^{th}$ derivative of $f_n$ with respect to $z$. From $(\ref{trigpol2})$
and $(\ref{zeros1})$ we get
\begin{eqnarray*}
0 & = & 1+a_{n,1} + a_{n,2} + a_{n,3} + \cdots + a_{n,M}\\
  & = & \lambda_n + a_{n,1}\cdot\Lambda_1 + a_{n,2}\cdot\Lambda_2 + a_{n,3}\cdot\Lambda_3 + \cdots + a_{n,M}\cdot\Lambda_M\\
  & = & \lambda_n^2 + a_{n,1}\cdot\Lambda_1^2 + a_{n,2}\cdot\Lambda_2^2 + a_{n,3}\cdot\Lambda_3^2 + \cdots + a_{n,M}\cdot\Lambda_M^2
\end{eqnarray*}
and so on, until the last equation ($0=f_n^{(M-1)}(0)$)
\[
0 =  \lambda_n^{M-1} + a_{n,1}\cdot\Lambda_1^{M-1} + a_{n,2}\cdot\Lambda_2^{M-1} + a_{n,3}\cdot\Lambda_3^{M-1} + \cdots +a_{n,M}\cdot\Lambda_M^{M-1}.
\]
We rewrite this system of equations as

\begin{equation}\label{matrixx}
\begin{pmatrix}
  -1 \\
  -\lambda_n  \\
  -\lambda_n^2 \\
 \vdots \\
 -\lambda_n^{M-1}
 \end{pmatrix}
 =
 \begin{pmatrix}
 1  & 1 & 1 &  \cdots & 1 \\
  \Lambda_1 & \Lambda_2 & \Lambda_3 & \cdots & \Lambda_M  \\
  \Lambda_1^2 & \Lambda_2^2 & \Lambda_3^2 & \cdots & \Lambda_M^2  \\
 \vdots  & \vdots  & \vdots & \ddots & \vdots \\

 \Lambda_1^{M-1} & \Lambda_2^{M-1} & \Lambda_3^{M-1} & \cdots & \Lambda_M^{M-1}
 \end{pmatrix}
 \cdot
\begin{pmatrix}
  a_{n,1} \\
  a_{n,2} \\
  a_{n,3}\\
  \vdots \\
  a_{n,M}
 \end{pmatrix}.
 \end{equation}
Since the $\Lambda_i$'s are all different, then the above square Vandermonde matrix is invertible. Therefore, real numbers
$a_{n,j}$, $j=1,2,\cdots, M$ do exist so that $f_n$ vanishes at $z=0$ at least $M$ times. Thus
\[
f_n(z)=z^M\cdot g_n(z)\qquad \text{for\,\, some\,\, entire\,\, function}\quad g_n(z),
\]
hence $(\ref{Mc})$ holds.

\subsection{The system $(\ref{system1})$ is Exact if and only if $M-\frac{1}{2}\le \alpha<M+\frac{1}{2}$}

If $\alpha\in [M-\frac{1}{2}, M+\frac{1}{2})$,
then the system $(\ref{system1})$ is both complete and minimal in $L^2 (0,1)$, hence it is exact.\\

\smallskip

Assume now that the system $(\ref{system1})$ is exact: we will show that $\alpha\in [M-\frac{1}{2}, M+\frac{1}{2})$.

Suppose first that it is exact but $\alpha <M-\frac{1}{2}$, say $M-\frac{3}{2}\le \alpha <M-\frac{1}{2}$.
Then, if $B\subset\mathbb{Z}$ such that $B=A\setminus\{\tau\}$ where $\tau$ is an integer not in the set $A$,
then the cardinality of $B$ equals $M-1$. Hence the system $\{t^{\alpha}\cdot r_n(t)\}_{n\in\mathbb{Z}\setminus B}$ is exact.
So the system $\{t^{\alpha}\cdot r_n(t)\}_{n\in\mathbb{Z}\setminus A}$
will have excess equal to $-1$ instead of zero. Similarly for other positive values of $\alpha$ when  $\alpha <M-\frac{3}{2}$.
So $\alpha$ cannot be less than $M-\frac{1}{2}$.

Suppose now that it is exact but $\alpha \ge M+\frac{1}{2}$, say $M+\frac{1}{2}\le \alpha <M+\frac{3}{2}$.
Then, if $C\subset\mathbb{Z}$ such that $C=A\cup\{\tau\}$ where $\tau$ is an integer not in the set $A$,
then the cardinality of $C$ equals $M+1$. Hence the system $\{t^{\alpha}\cdot r_n(t)\}_{n\in\mathbb{Z}\setminus C}$ is exact.
So the system $\{t^{\alpha}\cdot r_n(t)\}_{n\in\mathbb{Z}\setminus A}$
will have excess equal to $+1$ instead of zero. Similarly for other positive values of $\alpha$ when  $\alpha\ge M+\frac{3}{2}$.
So $\alpha$ cannot be greater than or equal to $M+\frac{1}{2}$.

In this way we conclude that if the system $(\ref{system1})$ is exact, then $\alpha\in [M-\frac{1}{2}, M+\frac{1}{2})$.

\subsection{The exact system $(\ref{system1})$ is not a Riesz Basis for $L^2 (0,1)$}

Suppose that for $\alpha\in [M-\frac{1}{2} , M+\frac{1}{2})$ the exact system
$\{t^{\alpha}\cdot r_n(t)\}_{n\in\mathbb{Z}\setminus A}$
is a Riesz Basis for $L^2 (0,1)$ and hence a Schauder Basis as well. Let $f_n$ be as in $(\ref{trigpol2})$:
since the family $\{f_n/t_{\alpha}\}_{n\in\mathbb{Z}\setminus A}$ is biorthogonal
to the Schauder Basis $(\ref{system1})$, then $\{f_n/t_{\alpha}\}_{n\in\mathbb{Z}\setminus A}$
is the dual Schauder basis, hence for any $f\in L^2 (0,1)$, we have
\[
f(t)=\sum_{n\in\mathbb{Z}\setminus A} \langle f , \frac{f_n}{t_{\alpha}}\rangle \cdot t^{\alpha} r_n(t),
\]
with convergence in the $L^2 (0,1)$ norm. This implies that $|\langle f , \frac{f_n}{t_{\alpha}} \rangle |\cdot || t_{\alpha} r_n||\to 0$ as $n\to\infty$.

In particular, for each integer $m$ belonging to the set $A$ which has been excluded from $\mathbb{Z}$,
we have
\begin{equation}\label{basisseries}
t^{\alpha}\cdot r_m(t)=\sum_{n\in\mathbb{Z}\setminus A} \langle t_{\alpha}r_m , \frac{f_n}{t_{\alpha}} \rangle \cdot t^{\alpha} r_n(t)
\end{equation}
with convergence in the $L^2 (0,1)$ norm. Fix such an $m\in A$. Then

\begin{equation}\label{tozero}
|\langle t_{\alpha}r_m , \frac{f_n}{t_{\alpha}} \rangle|\cdot || t_{\alpha}r_n||\to 0\qquad \text{as}\qquad n\to\infty.
\end{equation}

In $(\ref{trigpol2})$ we have for every $n\in\mathbb{Z}\setminus A$,
$f_n(t)=e^{i\lambda_n t}+a_{n,1}e^{i\Lambda_1 t}+a_{n,2}e^{i\Lambda_2 t}+\dots +a_{n,M}e^{i\Lambda_M t}$.
Since $m\in A$ and $\{\lambda_i:\,\, i\in A\}$ is rewritten as $\{\Lambda_1, \Lambda_2, \dots ,\Lambda_M\}$
such that $\Lambda_1<\Lambda_2<\dots <\Lambda_M$, then $\lambda_m=\Lambda_j$ for some $j\in\{1,2,\cdots, M\}$, and
$e^{i\lambda_m t}=e^{i\Lambda_j t}$.
Since the family $\{r_n\}_{n\in\mathbb{Z}}$ is biorthogonal to the family $\{e^{i\lambda_n t}\}_{n\in\mathbb{Z}}$, we get
\[
\langle t_{\alpha}r_m, \frac{f_n}{t_{\alpha}}\rangle= \overline{a_{n,j}}\qquad \text{for\,\, all}\quad n\in\mathbb{Z}\setminus A.
\]
Replacing into $(\ref{tozero})$ we have
\begin{equation}\label{tozero1}
|a_{n,j}|\cdot || t_{\alpha}r_n||\to 0\qquad \text{as}\qquad n\to\infty.
\end{equation}

Now, if the system $(\ref{system1})$ is a Riesz basis for $L^2 (0,1)$, then there are positive constants $p$ and $q$ so that
$p<|| t_{\alpha}r_n||<q$ for all $n\in\mathbb{Z}\setminus A$ (\cite[Lemma 3.6.9]{Christensen}).
We will show that there is a positive number $\delta$, independent of $n\in\mathbb{Z}\setminus A$ and $j\in\{1,2,\dots, M\}$, so that
\begin{equation}\label{delta}
\delta\cdot |\lambda_n|^{M-1}\le |a_{n,j}|\qquad \text{finally\,\, for\,\, all}\quad n\in\mathbb{Z}\setminus A,\quad \text{and} \quad j\in\{1,2,\dots, M\}.
\end{equation}
But then, combining the above  with $p<|| t_{\alpha}r_n||<q$ for all $n\in\mathbb{Z}\setminus A$, gives
\[
p\cdot\delta\cdot |\lambda_n|^{M-1}\le |a_{n,j}|\cdot || t_{\alpha}r_n||
\qquad \text{for\,\, all}\,\, n\in\mathbb{Z}\setminus A,\quad \text{and} \quad j\in\{1,2,\dots, M\}.
\]
However, $M-1\ge 0$ and $|\lambda_n|\to \infty$ as $n\to\infty$, thus the above shows that $(\ref{tozero1})$ is not true.
Hence the exact system $(\ref{system1})$ is not a Riesz Basis for $L^2 (0,1)$.

\smallskip

So it remains to verify $(\ref{delta})$. Returning to the relation $(\ref{matrixx})$, consider the invertible Vandermonde matrix
\[
\Lambda:=\begin{pmatrix}
 1  & 1 & 1 &  \cdots & 1 \\
  \Lambda_1 & \Lambda_2 & \Lambda_3 & \cdots & \Lambda_M  \\
  \Lambda_1^2 & \Lambda_2^2 & \Lambda_3^2 & \cdots & \Lambda_M^2  \\
 \vdots  & \vdots  & \vdots & \ddots & \vdots \\

 \Lambda_1^{M-1} & \Lambda_2^{M-1} & \Lambda_3^{M-1} & \cdots & \Lambda_M^{M-1}
 \end{pmatrix}
\]
and denote by $\Lambda^{-1}$, $det (\Lambda)$, and $adj (\Lambda)$, the inverse, determinant and adjoint of $\Lambda$ respectively.
Then we have
\[
\begin{pmatrix}
  -1 \\
  -\lambda_n  \\
  -\lambda_n^2 \\
 \vdots \\
 -\lambda_n^{M-1}
 \end{pmatrix}
 =\Lambda
 \cdot
\begin{pmatrix}
  a_{n,1} \\
  a_{n,2} \\
  a_{n,3}\\
  \vdots \\
  a_{n,M}
 \end{pmatrix}
 \]
and
\begin{eqnarray*}
\begin{pmatrix}
  a_{n,1} \\
  a_{n,2} \\
  a_{n,3}\\
  \vdots \\
  a_{n,M}
 \end{pmatrix}& = & \Lambda^{-1}\cdot \begin{pmatrix}
  -1 \\
  -\lambda_n  \\
  -\lambda_n^2 \\
 \vdots \\
 -\lambda_n^{M-1}
 \end{pmatrix}\\
& = &
 \frac{1}{det (\Lambda)}\cdot adj (\Lambda)\cdot
 \begin{pmatrix}
  -1 \\
  -\lambda_n  \\
  -\lambda_n^2 \\
 \vdots \\
 -\lambda_n^{M-1}
 \end{pmatrix}.
\end{eqnarray*}
We now claim that all the entries of the far right column of $adj (\Lambda)$ are non-zero real numbers.
So we need to check for the cofactors of the entries of the bottom row of $\Lambda$. For example, for the entry
$\Lambda_2^{M-1}$, its cofactor is equal to $(-1)^{2+M}$ times the determinant of the matrix that remains
when we delete the second column and the last row. That matrix is
\[
\begin{pmatrix}
 1  & 1 & 1 &  \cdots & 1 \\
  \Lambda_1 & \Lambda_3 &  \Lambda_4 & \cdots & \Lambda_M  \\
  \Lambda_1^2  & \Lambda_3^2 &  \Lambda_4^2 & \cdots & \Lambda_M^2  \\
 \vdots  & \vdots & \ddots & \vdots \\

 \Lambda_1^{M-2} & \Lambda_3^{M-2} &  \Lambda_4^{M-2} & \cdots & \Lambda_M^{M-2}
 \end{pmatrix}.
\]
Clearly this is again an invertible Vandermonde matrix. The same is true for the cofactors of the other elements of the last row of $\Lambda$,
thus we conclude that all the entries of the far right column of $adj (\Lambda)$ are non-zero real numbers.
This means that for all $n\in\mathbb{Z}\setminus A$ and all $j=1,2,\dots, M$, we have
\[
a_{n,j}=\sum_{k=1}^{M-1} c_{j,k} (-\lambda_n)^{k-1} + c_{j,M} (-\lambda_n)^{M-1}
\]
where all the $c_{j,k}$ are real numbers, and especially $c_{j,M}\not= 0$ for all $j=1,2,\dots, M$.
This implies the existence of some positive $\delta$ so that $(\ref{delta})$ holds.

\subsection{The exact system $(\ref{system})$ is not a Schauder Basis for $L^2 (0,1)$}

If the exponential Schauder basis for $L^2 (0,1)$ is the family $\{e^{i2\pi nt}\}_{n\mathbb{Z}}$,
then its dual Schauder basis is the same family, thus $r_n=e^{i2\pi nt}$ for all $n\in\mathbb{Z}$.
In this case, we have $|| t_{\alpha} r_n||=1$, thus $(\ref{tozero1})$ becomes
\[
|a_{n,j}|\to 0\qquad \text{as}\qquad n\to\infty.
\]
Clearly this is false due to $(\ref{delta})$, therefore the exact system $(\ref{system})$ is not a Schauder Basis for $L^2 (0,1)$.

The proof of Theorem $\ref{theorem}$ is now finished.

\end{document}